\documentclass[12pt]{article}

\usepackage{amsmath}
\usepackage{amsfonts}

\newtheorem{theorem}{Theorem}
\newtheorem{lemma}[theorem]{Lemma}
\newtheorem{definition}[theorem]{Definition}

\begin{document}

\title{Multidimensional $p$-adic wavelets \\ for the deformed metric}

\author{S.Albeverio, S.V.Kozyrev}

\maketitle

\begin{abstract}
The approach to $p$-adic wavelet theory from the point of view of representation theory is discussed. $p$-Adic wavelet frames can be constructed as orbits of some $p$-adic groups of transformations. These groups are automorphisms of the tree of balls in the $p$-adic space.
In the present paper we consider deformations of the standard $p$-adic metric in many dimensions and construct some corresponding groups of transformations.  We build several examples of $p$-adic wavelet bases. We show that the constructed wavelets are eigenvectors of some pseudodifferential operators.
\end{abstract}

Keywords:  $p$-adic wavelets, multiresolution
analysis, representation theory, spectral theory

\section{Introduction}

$p$-Adic wavelets were defined in \cite{wavelets} where also the relation between wavelets and spectral analysis of $p$-adic pseudodifferential operators was described. In \cite{ACHA} these results were generalized to wavelets and operators on general locally compact ultrametric spaces. In \cite{KSS}, \cite{KSS1}, \cite{A}, \cite{B}, \cite{C}, \cite{D}  $p$-adic multiresolution wavelet analysis was investigated.

Approach in $p$-adic wavelet analysis related to the representation theory of some $p$-adic groups of transformations was proposed in \cite{frames}, \cite{framesdimd}. All transformations from these group map balls to balls. $p$-Adic wavelet bases and frames in this approach are constructed as orbits of the groups of transformations.

In the multidimensional case the metric in $p$-adic spaces can be defined non-uniquely. Moreover the different metrics will correspond to the different trees of balls in $\mathbb{Q}_p^d$ with the different groups of automorphisms (automorphisms of the tree of balls in an ultrametric space are also called ball--morphisms). These metrics will be related to the different wavelet bases.
In the present paper we discuss some examples of metrics in $\mathbb{Q}_p^d$ and construct the corresponding wavelet bases. These wavelet bases are related to matrix dilations.

Matrix dilations were used in the real wavelet analysis to construct multiresolution wavelet bases in many dimensions, some of these bases contain wavelets supported on fractals \cite{GroMad}. The example of the corresponding matrix dilation is given by the quincunx matrix.

In paper \cite{KingSkopina} the example of the 2-adic 2-dimensional wavelet basis was constructed, which is related to the action of the quincunx matrix. This action does not map balls to balls, if we consider the standard 2-dimensional ultrametric. In the present paper we propose the example of metric in $\mathbb{Q}_2^2$ such that the quincunx matrix acts as a ball--morphism for this metric (i.e. maps balls to balls).

We also construct the family of general matrix dilations related to some non-standard metrics in $\mathbb{Q}_p^d$ and consider the corresponding wavelet bases. These bases contain the functions
$$
\Psi_{k;jn}(x)=p^{j\over 2}\Psi_{k}(A^{-j}x-n),\qquad j\in \mathbb{Z},\quad n\in \mathbb{Q}_p^d/\mathbb{Z}_p^d,
$$
$$
\Psi_k(x)=\chi\left(k\cdot A^{-1}x\right)\Omega(|x|_p),\quad k\in \mathbb{Z}_p^d/A^{*}\mathbb{Z}_p^d\backslash\{0\}, \quad k\cdot x=\sum_{i=1}^{d}k_ix_i,
$$
$A$ is the matrix dilation, see section 6 for the details. The $n$ and $k$ are understood as some (fixed in explicit way, see below) representatives from the corresponding equivalence classes. The above wavelets are analogues of the standard one dimensional $p$-adic wavelet $\chi(p^{-1}x)\Omega(|x|_p)$ introduced in \cite{wavelets}.

Moreover, the above wavelets are eigenfunctions of the natural pseudodifferential operators
$$
D^{\alpha}\Psi_{k;jn}=\|{A^*}^{(-j-1)}k\|^{\alpha}\Psi_{k;jn},
$$
where
$$
D^{\alpha}f(x)=F^{-1} \left(\|k\|^{\alpha}F[f]\right)(x),
$$
where $\|\cdot\|$ is the norm in $\mathbb{Q}_p^d$ related to the considered metric, $F$ is the Fourier transform.

Groups of automorphisms of trees were discussed in particular in \cite{Olshansky}, \cite{Neretin}, \cite{Serre}, \cite{Serre1}, \cite{Cartier1}. If the group conserves the path to infinity in the tree with the above partial order, the group is called parabolic, if the group conserves all elements of this path to infinity (starting from some vertex in the path) then the group is called orispheric.

The structure of the present paper is as follows.

In section 2 we discuss deformation of metric in multidimensional $p$-adic spaces and describe the corresponding groups of isometries.

In section 3 we discuss and construct some examples of dilations for the discussed in section 2 deformed metrics. A dilation here is a linear map which maps any ball with the center in zero to the maximal subball of this ball.

In sections 4 and 5 we consider two examples of deformed metrics, dilations and the corresponding wavelet bases in $L^2(\mathbb{Q}_2^2)$ (i.e. discuss the 2-dimensional 2-adic case). The example of section 5 is related to the quincunx matrix which in the real case generate wavelets supported on  fractals.

In section 6 we introduce the wavelet basis related to deformed metric and dilation in $\mathbb{Q}_p^d$ and show that introduced wavelets are eigenvectors of some pseudodifferential operators.

In section 7 we discuss the map of introduced $p$-adic wavelet to real spaces and show that the condition of paper \cite{GroMad} of existence of the corresponding real wavelet basis is equivalent to the statement that the corresponding natural map of $\mathbb{Q}_p^d$ to $\mathbb{R}^d$ conserves the measure (the measure of a set is equal to the measure of the image of this set).

\section{Deformed metrics and groups of isometries}

In the present section we consider deformations of the ultrametric in $\mathbb{Q}_p^d$ and find the corresponding groups of isometries. We will show that deformation of the metric changes the group of isometries.

The standard ultrametric on $\mathbb{Q}_p^d$ is defined by the formula
\begin{equation}\label{standard_norm}
d(x,y)=|x-y|_p={\rm max}(|x_l-y_l|_p),\quad l=1,\dots,d,
\end{equation}
$$
x=(x_1,\dots,x_d),\quad y=(y_1,\dots,y_d).
$$

Let us deform the above definition by an introduction of weights $q_l$
\begin{equation}\label{deformed_norm}
\widetilde{d}(x,y)={\rm max}(q_l|x_l-y_l|_p),\quad l=1,\dots,d,\quad p^{-1}<q_l\le 1.
\end{equation}
This defines a deformed ultrametric on $\mathbb{Q}_p^d$. The corresponding norm on $\mathbb{Q}_p^d$ we denote $\|\cdot\|$. We will use also more general deformed ultrametrics (which we will also denote $\widetilde{d}$) which can be obtained from (\ref{deformed_norm}) by linear transformations from the group $O_d$ of norm conserving linear transformations, see below. In particular, for any deformed metric the set $\{p^{j}\mathbb{Z}_p^d\}$, $j\in\mathbb{Z}$ is a subset of the set of balls with the center in zero.

The translations of $p^{j}\mathbb{Z}_p^d$, $j\in\mathbb{Z}$ (where $\mathbb{Z}_p^d$ is the unit ball with respect to (\ref{standard_norm})) give the the set ${\cal T}(\mathbb{Q}_p^d,d)$ of all balls for the metric (\ref{standard_norm}). These sets will be also balls with respect to any of the deformed ultrametrics (\ref{deformed_norm}). For deformed ultrametrics we will have also other examples of balls which will be the finite unions of maximal subballs in $p^{j}\mathbb{Z}_p^d$, see below.

Let us discuss the group of linear isometries for the metric (\ref{standard_norm}).

\begin{lemma}\quad{\sl
The group $O_d$ of linear transformations in $\mathbb{Q}_p^d$ which conserve the metric $d(\cdot,\cdot)$ given by (\ref{standard_norm}) coincides with the set of matrices with matrix elements in $\mathbb{Z}_p$ and $|{\rm det}(\cdot)|_p=1$.}
\end{lemma}

\noindent{\it Proof}\qquad
The volume $\mu({\cal W})$ of the $\mathbb{Z}_p$-module ${\cal W}$ generated by the set of vectors $x_i\in\mathbb{Q}_p^d$, $i=1,\dots,d$ is equal to $|{\rm det}(X)|_p$, where $X$ is the $d\times d$ matrix with the columns $x_i$. Here $\mu$ the Haar measure on $\mathbb{Q}_p^d$, the measure of the unit ball $\mathbb{Z}_p^d$ is equal to one.

When the matrix $X$ is diagonal this statement is obvious. In the general case, both the volume of the module ${\cal W}$ and $|{\rm det}(X)|_p$ are invariant with respect to the transpositions of the columns and to the adding to a column of a linear combination of other columns. Performing these transformations we can bring the matrix $X$ to the diagonal form, where $\mu({\cal W})=|{\rm det}(X)|_p$.

In \cite{framesdimd} it was proven that the group $O_d$ of norm conserving linear transformations coincides with the stabilizer (in the group of non degenerate linear transformations) of the $\mathbb{Z}_p$-module ${\cal V}$ generated by vectors from the coordinate basis in $\mathbb{Q}_p^d$ (i.e. by the columns with one element equal to one and other elements equal to zero). Therefore $O_d$ is a subgroup in the group of matrices with $|{\rm det}(\cdot)|_p=1$.

A matrix with elements in $\mathbb{Z}_p$ maps ${\cal V}$ to itself. Since the matrix with $|{\rm det}(\cdot)|_p=1$ conserves the volume $\mu({\cal V})$ this implies that the matrix with elements in $\mathbb{Z}_p$ and with $|{\rm det}(\cdot)|_p=1$ belongs to $O_d$.

It is easy the check that a matrix with any of the matrix elements belonging to $\mathbb{Q}_p\backslash \mathbb{Z}_p$ does not map ${\cal V}$ to itself. This finishes the proof of the lemma. $\square$

\bigskip

We make the following observation: the deformed metric (\ref{deformed_norm}) possesses different group of isometries in comparison to the standard metric  (\ref{standard_norm}). Moreover, for the deformed metric the group of isometries will be smaller --- it will be a subgroup of the group of isometries for the standard metric.

Making transpositions of the coordinates one can bring the deformed metric (\ref{deformed_norm}) in the form with $q_1\le q_2\le \dots \le q_d$. We have the following lemma which describes the group of isometries for the deformed metric.

\begin{lemma}\label{Borel}\quad
{\sl Let us consider the deformed metric (\ref{deformed_norm}) with the parameters $q_1< q_2< \dots < q_d$ (all the parameters $q_i$ belong to $(p^{-1},1]$).

The group of isometries $O(\widetilde{d})$ for this metric is a subgroup of $O_d$ which contains the matrices $A\in O_d$
which are
${\rm mod}\,p$ upper triangular, i.e. the matrix elements $A_{ab}$ with $a>b$ belong to $p\mathbb{Z}_p$.}
\end{lemma}

\noindent{\it Proof}\qquad  The sequence of balls between $p\mathbb{Z}_p^d$ and $\mathbb{Z}_p^d$ for the described metric $\widetilde{d}$ is the set of products
\begin{equation}\label{sequence_balls}
B_a=\mathbb{Z}_p\times\dots\mathbb{Z}_p\times p\mathbb{Z}_p\dots p \mathbb{Z}_p
\end{equation}
with $a$ components $\mathbb{Z}_p$ and $d-a$ components $p\mathbb{Z}_p$, $a=0,\dots, d$.

All the above subballs (containing zero) have different diameters, namely the diameter of the subball with the parameter $a$ is equal to $q_a$ for $a=1,\dots,d$ and is equal to $p^{-1}q_d$ for $a=0$. All the balls in $\mathbb{Q}_p^d$ (with respect to the mentioned above deformed metric) are translations and dilations by degrees of $p$ of the described set of balls.

An isometry from $O(\widetilde{d})$ should conserve the above sequence of balls. The group of isometries in an ultrametric space coincides with the group of orispheric automorphisms of the corresponding tree of balls in this space (where orispheric map is the map which conserves all the balls in some maximal increasing sequence of balls $I_0<I_1<I_2<\dots$ where $I_0$ depends on the map), see for example \cite{vectorPDO}. In particular, transformations from $O(\widetilde{d})$ act by orispheric automorphisms of ${\cal T}(\mathbb{Q}_p^d,\widetilde{d})$.
Therefore $O(\widetilde{d})$ is a subgroup of $O_d$ (since $O(\widetilde{d})$ conserves the balls $p^j\mathbb{Z}_p^d$).

Moreover $O(\widetilde{d})$ is the stabilizer of the set of balls (\ref{sequence_balls}) in $O_d$ (since any map from the stabilizer is a linear orispheric automorphism of the tree ${\cal T}(\mathbb{Q}_p^d,\widetilde{d})$ of balls and therefore belongs to $O(\widetilde{d})$).

Since $\mathbb{Z}_p/p\mathbb{Z}_p=\mathbb{F}_p$ is the finite field of residues ${\rm mod}\, p$, the above sequence of balls (\ref{sequence_balls}) can be considered as a complete flag in $\mathbb{F}_p^d$, i.e. a maximal increasing sequence of linear subspaces of the $d$-dimensional linear space over the finite field $\mathbb{F}_p$. Namely factorizing a ball from (\ref{sequence_balls}) over $p\mathbb{Z}_p^d$ we get a linear subspace in $\mathbb{F}_p^d$.

The stabilizer of the  complete flag in $\mathbb{F}_p^d$ is the group of invertible upper triangular matrices over $\mathbb{F}_p$ (a Borel subgroup in the group of non degenerate matrices over $\mathbb{F}_p$). This group is the image of $O(\widetilde{d})$ with respect to factorization ${\rm mod}\, p\mathbb{Z}_p$.

Therefore $O(\widetilde{d})$ coincides with the subgroup of $O_d$ of ${\rm mod}\,p$ upper triangular matrices.
This finishes the proof of the lemma. $\square$

\bigskip

The next lemma describes isometries for a more general deformed ultrametric.

\begin{lemma}\quad
{\sl
Let us consider the deformed metric (\ref{deformed_norm}) with the parameters $q_1\le q_2\le \dots \le q_d$ from $(p^{-1},1]$, and, moreover, the coinciding parameters $q_a$ look as follows: $q_{a}=q_{b}$ for $d_{i-1}<a,b\le d_{i}$, where $d_i$ are some integers belonging to the set $\{1,\dots,d\}$.

The group of isometries $O(\widetilde{d})$ for this metric is a subgroup of $O_d$ which contains the matrices $A\in O_d$ which are block ${\rm mod}\,p$ upper triangular, i.e. the matrix elements $A_{ab}$ with $a>b$ and $a$, $b$ which do not lie in some $(d_{i-1},d_i]$ (i.e. the corresponding $q_a$ and $q_b$ are different) belong to $p\mathbb{Z}_p$.
}
\end{lemma}

The proof of this lemma is a direct generalization of the proof of the previous lemma to the case of incomplete flags. We recall that an incomplete flag in $\mathbb{F}_p^d$ is a not necessarily a maximal sequence of increasing linear subspaces in $\mathbb{F}_p^d$ which contains $\{0\}$ and $\mathbb{F}_p^d$.

\section{Dilations}

To construct multidimensional wavelet bases we use the following definition of dilations. Let $\widetilde{d}$ be an arbitrary deformed ultrametric in $\mathbb{Q}_p^d$.

\begin{definition}\label{dilation}\quad{\sl
A dilation with respect to a deformed ultrametric $\widetilde{d}$ is a linear map $\mathbb{Q}_p^d\to \mathbb{Q}_p^d$ which maps any $\widetilde{d}$-ball with the center in zero to a maximal $\widetilde{d}$-subball (with the center in zero) of this ball.
}
\end{definition}

Let us note that for a dilation the above condition must hold for {\it any} ball with the center in zero. Therefore for a dilation $A$ the set of balls with the center in zero is given by $A^{j}\mathbb{Z}_p^d$, $j\in \mathbb{Z}$ (since $\mathbb{Z}_p^d$ is a ball).

Since any ball is translation of a ball with the center in zero it follows that a dilation is a ball--morphism (an automorphism of the tree of balls ${\cal T}(\mathbb{Q}_p^d,\widetilde{d})$).

In particular, for the case of the deformed ultrametric with the parameters $q_1< q_2< \dots < q_d$ (which corresponds to a complete flag in the above Lemma \ref{Borel} the dilation should satisfy
$|{\rm det}(\cdot)|_p=p^{-1}$.

\begin{lemma}\label{orbit}\quad{\sl Let $A$ be a dilation in $\mathbb{Q}_p^d$ with respect to a deformed metric $\widetilde{d}$.

Then the set of characteristic functions of balls with respect to the deformed metric $\widetilde{d}$ is in one to one correspondence with the set of functions
$$
\Omega(|A^{j}x-n|_p)=\Omega(\|A^{j}x-n\|),\qquad j\in \mathbb{Z},\quad n\in \mathbb{Q}_p^d/\mathbb{Z}_p^d.
$$
Here $\Omega(\cdot)$ is the characteristic function of the unit interval $[0,1]$, the group $\mathbb{Q}_p^d/\mathbb{Z}_p^d$  is understood as the group of d-tuples of fractions with the addition modulo one:
$$
n=(n_1,\dots, n_d),\qquad n_i=\sum_{k=\beta_i}^{-1}n^{(i)}_k p^k,\quad i=1,\dots,d.
$$

}
\end{lemma}

\noindent{\it Proof}\qquad
The group of degrees of the dilation $A$ is transitive on the set of $\widetilde{d}$-balls which contain zero. The translations from $\mathbb{Q}_p^d/\mathbb{Z}_p^d$ applied to $\mathbb{Z}_p^d$ give all the $\widetilde{d}$-balls with the diameter 1.
This proves the lemma.
$\square$

\bigskip

Let us construct the example of a dilation in $\mathbb{Q}_p^d$. Consider the product of $d\times d$ matrices
$$
A=\left(\begin{array}{ccccc}
                           0 & 1 & 0 & \dots & 0\cr
                           0 & 0 & 1 & \dots & 0\cr
                           \vdots & & & & \vdots\cr
                           0 & 0 & \dots & 0 & 1\cr
                           p & 0 & \dots & 0 & 0\cr
                           \end{array}\right)=
\left(\begin{array}{cccc}1 & 0 & \dots & 0\cr
                           0 & 1 & \dots & 0\cr
                           \vdots & & & \vdots\cr
                           0 & \dots & 1 & 0 \cr
                           0 & \dots & 0 & p \cr
                           \end{array}\right)
                           \left(\begin{array}{ccccc}
                           0 & 1 & 0 & \dots & 0\cr
                           0 & 0 & 1 & \dots & 0\cr
                           \vdots & & & & \vdots\cr
                           0 & 0 & \dots & 0 & 1\cr
                           1 & 0 & \dots & 0 & 0\cr
                           \end{array}\right).
$$
The second matrix in the above product is the cyclic substitution of the standard coordinate basis.

\begin{lemma}\label{Borel_dilation}\quad{\sl The above matrix $A$ is a dilation for $\mathbb{Q}_p^d$
with the metric $\widetilde{d}$ described in lemma \ref{Borel}.
}
\end{lemma}

\noindent{\it Proof}\qquad
Let us consider the action of $A$ on the balls $B_a$ given by (\ref{sequence_balls}). We get
$$
AB_a=B_{a-1},\qquad a=1,\dots,d.
$$

Since dilations by $p^{j}$, $j\in\mathbb{Z}$ of $B_a$, $a=1,\dots,d$, constitute all $\widetilde{d}$-balls containing zero this proves the lemma. $\square$

\bigskip

One can use dilations (together with translations and the isometries described in the present section) to construct wavelet bases as orbits of groups generated by the mentioned transformations as in \cite{framesdimd}, \cite{frames}.
In the next two sections we will consider the examples of dilations and wavelet bases discussed in \cite{GroMad} in the framework of the real wavelet analysis (the first example actually is the particular case of the above lemma).

\section{Example of the two dimensional wavelet basis}

In the present section we consider the example of the wavelet basis in $L^2(\mathbb{Q}_2^2)$.

Consider the matrix
\begin{equation}\label{S}
S=\left(\begin{array}{cc}0 & 1\cr 2 & 0\cr\end{array}\right).
\end{equation}
It is easy to check that $S^2=2E$, where $E$ is the unit matrix, and ${\rm det}\, S=-2$.

The map $S$ is not a dilation of $\mathbb{Q}_2^2$ with respect to the standard metric, but this is a dilation with respect to the deformed metric, see below.

Let us consider for $\mathbb{Q}_2^2$ the deformed ultrametric with $q_2=1$, $q_1=q$, $p^{-1}<q<1$. This particular choice of the deformed metric on $\mathbb{Q}_2^2$ will be denoted by $s$.

\begin{lemma}\label{ballmorphism}\quad {\sl
The matrix $S$ is a dilation for the deformed ultrametric $s$.
}
\end{lemma}

\noindent{\it Proof}\qquad The set of all balls in $\mathbb{Q}_2^2$ with respect to the metric $s$ (or $s$-balls) contains the sets $2^{-a}\mathbb{Z}_2^2$, $a\in\mathbb{Z}$ and translations of these sets, and the sets
$$
2^{-b}\left(\begin{array}{c}\mathbb{Z}_2\cr 2\mathbb{Z}_2 \cr\end{array}\right),\quad b\in\mathbb{Z},
$$
and translations of these sets.

Here
$$
\left(\begin{array}{c}\mathbb{Z}_2\cr 2\mathbb{Z}_2 \cr\end{array}\right)
$$
is the set of $(x_1,x_2)\in\mathbb{Q}_2^2$ with $x_1\in \mathbb{Z}_2$ and $x_2\in 2\mathbb{Z}_2$.

The action of the matrix $S$ on balls in $\mathbb{Q}_2^2$ has the form
\begin{equation}\label{S1}
S\mathbb{Z}_2^2=\left(\begin{array}{c}\mathbb{Z}_2\cr 2\mathbb{Z}_2 \cr\end{array}\right)  ,\quad
S\left(\begin{array}{c}\mathbb{Z}_2\cr 2\mathbb{Z}_2 \cr\end{array}\right)=
\left(\begin{array}{c}2\mathbb{Z}_2\cr 2\mathbb{Z}_2 \cr\end{array}\right)=2\mathbb{Z}^2_2,
\end{equation}
\begin{equation}\label{S2}
S\left(\begin{array}{c}\mathbb{Z}_2\cr 1+2\mathbb{Z}_2 \cr\end{array}\right)=
\left(\begin{array}{c}1+2\mathbb{Z}_2\cr 2\mathbb{Z}_2 \cr\end{array}\right)=
\left(\begin{array}{c}1\cr 0 \cr\end{array}\right)+2\mathbb{Z}^2_2.
\end{equation}

Formulas (\ref{S1}), (\ref{S2}) imply that the application of $S$ to any of the balls belongs to the described family of balls. This finishes the proof of the lemma. $\square$

\bigskip

It is possible to describe all dilations for the metric $s$.

\begin{lemma}\label{s-dilations}\quad
{\sl The matrix
$$
A=\left(\begin{array}{cc}a & b\cr c & d\cr\end{array}\right)
$$
is a dilation for the metric $s$ if and only if
\begin{equation}\label{parameters}
a=0~\,{\rm mod}\, 2,\quad b=1~\,{\rm mod}\, 2,\quad c=2~\,{\rm mod}\, 4,\quad d=0~\,{\rm mod}\, 2.
\end{equation}
}
\end{lemma}

\noindent{\it Proof}\qquad The proof of this lemma is by consideration of the action of $A$ on the balls
$$
\mathbb{Z}_2^2,\qquad \left(\begin{array}{c}\mathbb{Z}_2\cr 2\mathbb{Z}_2 \cr\end{array}\right).
$$

Since
$$
A:\mathbb{Z}_2^2\mapsto \left(\begin{array}{c}\mathbb{Z}_2\cr 2\mathbb{Z}_2 \cr\end{array}\right)
$$
we get
$$
a~\,{\rm mod}\, 2\ne b~\,{\rm mod}\, 2,\qquad c~\,{\rm mod}\, 2= d~\,{\rm mod}\, 2.
$$
Then we take into account
$$
A:\left(\begin{array}{c}\mathbb{Z}_2\cr 2\mathbb{Z}_2 \cr\end{array}\right)\mapsto 2\mathbb{Z}_2^2
$$
which implies
$$
a=0~\,{\rm mod}\, 2,\quad b=1~\,{\rm mod}\, 2,\quad c=0~\,{\rm mod}\, 2,\quad d=0~\,{\rm mod}\, 2.
$$

Taking into account $|{\rm det}\, A|_2=1/2$ we get the condition $c=2~\,{\rm mod}\, 4$.

We have found that the set of conditions (\ref{parameters}) is necessary. It is easy to check that these conditions are sufficient for $A$ to be a $s$-dilation.

This finishes the proof of the lemma. $\square$

\bigskip

Let us consider the wavelet function of the form of the difference of the two characteristic functions of $s$-balls:
\begin{equation}\label{theta}
\theta(x)=\Omega(|S^{-1}x|_2)-\Omega\left(\left|S^{-1}x-\left(\begin{array}{c}1/2\cr 0\cr\end{array}\right)\right|_2\right).
\end{equation}

The following theorem describes the corresponding wavelet basis.

\begin{theorem}\label{wavelets_s}\quad{\sl
The set of functions $\theta_{jn}(x)=p^{j\over 2}\theta(S^{-j}x-n)$, $j\in \mathbb{Z}$, $n\in \mathbb{Q}_2^2/\mathbb{Z}_2^2$ is an orthonormal basis in $L^2(\mathbb{Q}_2^2)$.}
\end{theorem}

\noindent{\it Proof}\qquad The proof of this theorem is straightforward. Orthogonality of the wavelets can be checked directly and in essence follows from (\ref{S1}), (\ref{S2}).  Lemma \ref{orbit} implies that to prove the completeness of the system of wavelets it it sufficient to check Parseval's identity for the characteristic function of the unit ball, which can be done directly as in \cite{wavelets}.

This finishes the proof of the theorem. $\square$

\section{The quincunx wavelet basis}

In the present section we discuss, using the language of ball--morphisms,  the quincunx 2-adic wavelet basis which was considered in \cite{KingSkopina}.

The quincunx matrix $Q$ with ${\rm det}\,Q=2$
\begin{equation}\label{quincunx}
Q=\left(\begin{array}{cc}1 & -1\cr 1 & 1\cr\end{array}\right)
\end{equation}
$$
Q^2=\left(\begin{array}{cc}0 & -2\cr 2 & 0\cr\end{array}\right),\quad Q^4=-4E.
$$
(where $E$ is the unit matrix) was used in the real wavelet analysis for the construction of wavelets with supports on fractals \cite{GroMad}.

In \cite{KingSkopina} it was found that the 2-adic 2-dimensional wavelet basis with dilations generated by the matrix $Q$ consists of mean zero test functions (test functions are linear combinations of characteristic functions of balls, therefore the fractal real basis of \cite{GroMad} corresponds to the regular $p$-adic basis).

The next lemma shows that the quincunx matrix is a dilation for the metric $q$  on $\mathbb{Q}_2^2$ which is a rotation of the deformed metric $s$ defined in the previous section  (i.e. this lemma is the analogue of lemma \ref{ballmorphism} for the matrix $Q$).

\begin{lemma}\quad{\sl
Let us consider the metric $q$ on $\mathbb{Q}_2^2$ given by application to the deformed ultrametric $s$ of the linear transformation $U\in O_d$ of the form
\begin{equation}\label{U}
U=\left(\begin{array}{cc}1 & 0\cr 1 & 1\cr\end{array}\right),
\end{equation}
namely the metric which has the form
\begin{equation}\label{metric_q}
q(x,y)=s(Ux,Uy)={\rm max}(q|x_1-y_1|_p,|x_1+x_2-y_1-y_2|_p),\qquad p^{-1}<q<1.
\end{equation}

Then the quincunx matrix $Q$ is a dilation of $\mathbb{Q}_2^2$ with the metric $q$.
}
\end{lemma}

\noindent{\it Proof}\qquad The set of balls with respect to the metric $q$ is given by translations and dilations by the degrees of 2 of the sets
$$
\mathbb{Z}_2^2,\qquad 2\mathbb{Z}_2^2 \bigcup \left(2\mathbb{Z}_2^2+\left(\begin{array}{c}1\cr 1\cr\end{array}\right)\right).
$$

Let us discuss the action of the quincunx matrix $Q$. The image of the unit ball $\mathbb{Z}_2^2$ with respect to $Q$ is the disjoint union of the two balls of the diameter $1/2$, and the image of $\mathbb{Z}_2^2$ w.r.t. $Q^2$ is $2\mathbb{Z}_2^2$:
\begin{equation}\label{Q}
Q:\mathbb{Z}_2^2 \to 2\mathbb{Z}_2^2 \bigcup \left(2\mathbb{Z}_2^2+\left(\begin{array}{c}1\cr 1\cr\end{array}\right)\right),
\end{equation}
\begin{equation}\label{Q1}
Q: 2\mathbb{Z}_2^2 \bigcup \left(2\mathbb{Z}_2^2+\left(\begin{array}{c}1\cr 1\cr\end{array}\right)\right)\to 2\mathbb{Z}_2^2,
\end{equation}
\begin{equation}\label{Q2}
Q: \left(2\mathbb{Z}_2^2+\left(\begin{array}{c}0\cr 1\cr\end{array}\right)\right) \bigcup \left(2\mathbb{Z}_2^2+\left(\begin{array}{c}1\cr 0\cr\end{array}\right)\right)\to 2\mathbb{Z}_2^2+\left(\begin{array}{c}1\cr 1\cr\end{array}\right),
\end{equation}
\begin{equation}\label{QQ}
Q^2:\mathbb{Z}_2^2 \to 2\mathbb{Z}_2^2.
\end{equation}

Taking into account formulas (\ref{Q}), (\ref{Q1}), (\ref{Q2}), (\ref{QQ}) we prove that $Q$ is a dilation with respect to the metric $q$.

This finishes the proof of the lemma. $\square$

\bigskip

\noindent{\bf Remark}\quad Since the metric $q$ given by (\ref{metric_q}) has the form $q(x,y)=s(Ux,Uy)$ and the matrix $S$ of the form (\ref{S}) is a dilation for the metric $s$, the matrix
$$
U^{-1}SU=\left(\begin{array}{cc}1 & 1\cr
                                1 & -1\cr\end{array}\right),
$$
where $U$ is given by (\ref{U}) is a dilation for the metric $q$.

Note that the quincunx matrix $Q$ does not coincide with the above example of a dilation for $q$ (but $UQU^{-1}$ fits in the conditions of lemma \ref{s-dilations} and therefore is a $s$-dilation).

\bigskip

Let us consider the wavelet equal to the difference of the two $q$--balls
\begin{equation}\label{psi}
\psi(x)=\Omega\left(|Q^{-1}x|_2\right)-\Omega\left(\left|Q^{-1}x-\left(\begin{array}{c}1/2\cr 1/2\cr\end{array}\right)\right|_2\right).
\end{equation}

The following theorem about the 2-adic quincunx wavelet basis has been proven in \cite{KingSkopina}.

\begin{theorem}\quad{\sl
The set of functions $\psi_{jn}(x)=p^{j\over 2}\psi(Q^{-j}x-n)$, $j\in \mathbb{Z}$, $n\in \mathbb{Q}_2^2/\mathbb{Z}_2^2$ is an orthonormal basis in $L^2(\mathbb{Q}_2^2)$.}
\end{theorem}

Alternatively this theorem can be proven in the same way as theorem \ref{wavelets_s} in the previous section.

\section{Wavelets and pseudodifferential operators}

Let $\widetilde{d}$ be some deformed ultrametric in $\mathbb{Q}_p^d$ and $A$ be a dilation with respect to $\widetilde{d}$, $|{\rm det}\, A|_p=p^{-1}$. In particular one can consider the metric described in lemma \ref{Borel} and the dilation constructed in lemma \ref{Borel_dilation}.

We define the wavelet function
\begin{equation}\label{wavelet1}
\Psi_k(x)=\chi\left(k\cdot A^{-1}x\right)\Omega(|x|_p),\quad k\in \mathbb{Z}_p^d/A^{*}\mathbb{Z}_p^d, \quad k\cdot x=\sum_{i=1}^{d}k_ix_i
\end{equation}
where $k$ does not represent zero in $\mathbb{Z}_p^d/A^{*}\mathbb{Z}_p^d$ (as usual $k$ are understood as some fixed representatives from the corresponding equivalence classes). Here $A^{*}$ is the transponated to $A$ matrix (for a dilation $A$ the transponated matrix is also a dilation).
The function $\chi$ is the additive character of $\mathbb{Q}_p$ (given by the expression (\ref{chi}) below).

We have $p-1$ wavelet functions of the above form (since there are $p$ representatives in $\mathbb{Z}_p^d/A^{*}\mathbb{Z}_p^d$). This definition is an analogue of the one dimensional $p$-adic wavelet $\chi(p^{-1}x)\Omega(|x|_p)$.

The parameter $k$ which enumerates the functions $\Psi_k$ can be considered as a parameter in the set of maximal subballs in the unit sphere $\|k\|=1$. If we investigate wavelets for the metric used in lemma \ref{Borel} then this set of $k$ can be considered as the parameter on the orbit of the group of isometries described in lemma \ref{Borel}.

The following formula is a generalization of (\ref{theta}), (\ref{psi}).

\begin{lemma}\label{wavelet_expand}\quad{\sl
Definition (\ref{wavelet1}) of the wavelet $\Psi_k$ can be equivalently rewritten as
\begin{equation}\label{expand}
\Psi_k(x)=\sum_{l=0}^{p-1}\chi\left(k\cdot A^{-1}m_l\right)\Omega\left(|A^{-1}(x-m_l)|_p\right),
\end{equation}
where $m_l$, $l=0,\dots, p-1$ are representatives in $\mathbb{Z}_p^d/A\mathbb{Z}_p^d$.

The coefficients $\chi\left(k\cdot A^{-1}m_l\right)$ are $p$--roots of one and  $\chi\left(k\cdot A^{-1}m_l\right)\ne \chi\left(k\cdot A^{-1}m_{l'}\right)$ for $l\ne l'$.

The wavelet $\Psi_k$ is a mean zero function.
}
\end{lemma}

\noindent{\it Proof}\qquad  Since
$$
\Psi_k(x)=\sum_{l=0}^{p-1}\chi\left(k\cdot A^{-1}(m_l+x-m_l)\right)\Omega\left(|A^{-1}(x-m_l)|_p\right)
$$
and in the above expression $A^{-1}(x-m_l)\in \mathbb{Z}_p^d$ (taking into account the support of $\Omega\left(|A^{-1}(\cdot-m_l)|_p\right)$) we get (\ref{expand}).

It obvious that
$$
\chi\left(k\cdot A^{-1}m_l\right)^p=\chi\left(pk\cdot A^{-1}m_l\right)=1.
$$

Since $|{\rm det}\, A|_p=p^{-1}$ the groups  $\mathbb{Z}_p^d/A\mathbb{Z}_p^d$ and $\mathbb{Z}_p^d/A^{*}\mathbb{Z}_p^d$ are isomorphic to the group of ${\rm mod}\, p$--residues. Therefore for $l\ne l'$ and $k\in \mathbb{Z}_p^d/A^{*}\mathbb{Z}_p^d\backslash{\{0\}}$
$$
\chi\left(k\cdot A^{-1}(m_l-m_{l'})\right)\ne 1.
$$

Since the function $\Psi_k$ takes values equal to the different $p$--roots of one on maximal $\widetilde{d}$--subballs in $\mathbb{Z}_p^d$ and the sum of all $p$--roots of one is zero this function has zero mean.

This finishes the proof of the lemma. $\square$

\bigskip

We construct wavelets using translations and dilations of $\Psi_k$:
\begin{equation}\label{wavelet2}
\Psi_{k;jn}(x)=p^{j\over 2}\Psi_{k}(A^{-j}x-n),\qquad j\in \mathbb{Z},\quad n\in \mathbb{Q}_p^d/\mathbb{Z}_p^d.
\end{equation}
As usual we understand the elements of the group $\mathbb{Q}_p^d/\mathbb{Z}_p^d$ as vectors with the coordinates belonging to the set of fractions as described in lemma \ref{orbit}.

\begin{theorem}\quad
{\sl The set of functions $\{\Psi_{k;jn}\}$ defined by (\ref{wavelet1}), (\ref{wavelet2}) is an orthonormal basis in $L^2(\mathbb{Q}_p^d)$.}
\end{theorem}

\noindent{\it Proof}\qquad  The proof in essence follows the proof of \cite{wavelets} and theorem \ref{wavelets_s} above.

Wavelets with the different indices $j$ (the scales of the wavelets) are orthogonal since one of the wavelets is constant on the support of the other wavelet and any wavelet is a mean zero function.

The supports of the wavelets with the same scale $j$ and the different indices $n$ (translations) have empty intersection.

Wavelets with the same $j$ and $n$ and the different $k$ are orthogonal by lemma \ref{wavelet_expand}.

The completeness of the system $\{\Psi_{k;jn}\}$ follows from the Parseval's identity. It is sufficient to prove this identity for the characteristic function $\Omega(|\cdot|_p)$ of the unit ball. One has
$$
\sum_{kjn}|\langle\Omega(|\cdot|_p) ,\Psi_{k;jn}\rangle|^2=\sum_{k;j<0}|\langle\Omega(|\cdot|_p) ,\Psi_{k;j0}\rangle|^2=(p-1)\sum_{j<0}p^{j}=1.
$$

This finishes the proof of the theorem. $\square$

\bigskip

Let us consider the pseudodifferential operator of the form
\begin{equation}\label{Dalpha}
D^{\alpha}f(x)=F^{-1} \left(\|k\|^{\alpha}F[f]\right)(x),
\end{equation}
where $\|\cdot\|$ is the $\widetilde{d}$-norm in $\mathbb{Q}_p^d$, $F$ is the Fourier transform
$$
F[f](k)=\int_{\mathbb{Q}_p^d} \chi(kx)f(x) d\mu(x),\quad k\cdot x=\sum_{i=1}^{d}k_ix_i,
$$
where $\mu$ is the ($d$-dimensional) Haar measure and $\chi$ is the additive character of $\mathbb{Q}_p$: $\chi(x)=\exp(2\pi i\{x\})$, $\{x\}$ is the fractional part of $x$:
\begin{equation}\label{chi}
x=\sum_{i=\gamma}^{\infty}x_ip^i,\quad x_i=0,\dots,p-1,\qquad \{x\}=\sum_{i=\gamma}^{-1}x_ip^i.
\end{equation}

This operator generalizes the Vladimirov fractional operator for the case of a deformed metric.

\begin{lemma}\quad{\sl
The basis $\{\Psi_{k;jn}\}$ is a basis of eigenvectors for the operator $D^{\alpha}$:
$$
D^{\alpha}\Psi_{k;jn}(x)=\|{A^*}^{(-j-1)}k\|^{\alpha}\Psi_{k;jn}(x),\quad k\in \mathbb{Z}_p^d/A^{*}\mathbb{Z}_p^d\backslash\{0\},\quad j\in \mathbb{Z},\quad n\in \mathbb{Q}_p^d/\mathbb{Z}_p^d.
$$
Here $A^*$ is the transponated to $A$ matrix.
}
\end{lemma}

\noindent{\it Proof}\qquad
Let us compute the Fourier transform of the wavelet functions:
$$
F[\Psi_{l}](k)=\int_{\mathbb{Q}_p^d} \chi(kx)\chi\left(l\cdot A^{-1}x)\right)\Omega(|x|_p) d\mu(x) =$$ $$=
\int_{\|x\|\le 1} \chi\left((k+{A^*}^{-1}l)\cdot x)\right) d\mu(x)=\Omega(|k+{A^*}^{-1}l|_p).
$$
Here $k\in \mathbb{Q}_p^d$.

Since
$$
\|k\|^{\alpha}\Omega(|k+{A^*}^{-1}l|_p)=\|{A^*}^{-1}l\|^{\alpha}\Omega(|k+{A^*}^{-1}l|_p)
$$
where $l \ne 0$ in $\mathbb{Z}_p^d/A^{*}\mathbb{Z}_p^d$,
we get
$$
D^{\alpha}\Psi_{l}=\|{A^*}^{-1}l\|^{\alpha}\Psi_{l}.
$$

Therefore applying $D^{\alpha}$ to wavelets (\ref{wavelet2}) we get
$$
D^{\alpha}\Psi_{l;jn}(x)=D^{\alpha} p^{j\over 2}\Psi_{l}(A^{-j}x-n)=\|{A^*}^{(-j-1)}l\|^{\alpha}\Psi_{l;jn}(x)
$$
since
$$
F[f(A\cdot)](k)=F[f]({A^*}^{-1}k)|{\rm det}\, A|_p^{-1}.
$$

This finishes the proof of the lemma. $\square$

\section{Map to real wavelets}

In \cite{wavelets} the following measure conserving 1--Lipschitz map (the Monna map, or the $p$-adic change of variable) of the $p$-adic field to the positive half--line was considered and it was proven that the images with respect to this map (for $p=2$) of 2-adic wavelets are the real Haar wavelets on the positive real half--line:
$$
\eta:\mathbb{Q}_p\to [0,\infty),
$$
$$
\eta: x=\sum_{i=\gamma}^{\infty}x_ip^i\mapsto \sum_{i=\gamma}^{\infty}x_ip^{-i-1}  ,\qquad x_i=0,\dots,p-1.
$$

In particular the above map maps $\mathbb{Z}_p$ to $[0,1]$ and $\mathbb{Q}_p/\mathbb{Z}_p$ (understood as the set of fractions $n=\sum_{\beta}^{-1}n_ip^i$, $n_i=0,\dots,p-1$ with the addition modulo one) to the set $\mathbb{Z}_{+}$ of non-negative integers. Note that the mentioned choice of the digits $n_i=0,\dots,p-1$ is important: if we choose some other set of digits then the set $\{n=\sum_{\beta}^{-1}n_ip^i\}$ will be the set of all representatives from $\mathbb{Q}_p/\mathbb{Z}_p$ but the image $\eta(\mathbb{Q}_p/\mathbb{Z}_p)$ of this set of representatives will be not equal to $\mathbb{Z}_{+}$. For example, with the choice of the digits $n_i=0,1+p,2(1+p),\dots,(p-1)(1+p)$ we will get $\eta(\mathbb{Q}_p/\mathbb{Z}_p)=(1+p)\mathbb{Z}_{+}$, $\eta(\mathbb{Z}_p)=[0,1+p]$.

In the present section we discuss the generalization of the above map to many dimensions and application to the construction of relation between multidimensional real and $p$-adic wavelets.

Consider the deformed metric $\widetilde{d}$ of $\mathbb{Q}_p^d$ (with $\mathbb{Z}_p^d$ being a ball) and the dilation $A\mathbb{Q}_p^d\to\mathbb{Q}_p^d$ with respect to the deformed metric.

A ball with the center in zero (with respect to the deformed metric) is an additive group.
Consider the factor group $\mathbb{Z}_p^d/A\mathbb{Z}_p^d$. This is a finite abelian group. We fix an arbitrary set of digits $\{n_0,\dots,n_{q-1}\}$, $q=|{\rm det}\, A|_p^{-1}$ --- vectors with integer coordinates which represent all elements of the factor group $\mathbb{Z}_p^d/A\mathbb{Z}_p^d$. The set $\{n_0,\dots,n_{q-1}\}$ of digits with the addition ${\rm mod}\, A\mathbb{Z}_p^d$ is isomorphic to the factor group  $\mathbb{Z}_p^d/A\mathbb{Z}_p^d$. For simplicity we consider the case $q=p$.

We have the following lemma.
\begin{lemma}\quad{\sl Let $A$ be a dilation in $\mathbb{Q}_p^d$ with respect to the metric $\widetilde {d}$ in the sense of definition \ref{dilation} and $|{\rm det}\, A|_p=p^{-1}$. Then the group $\mathbb{Z}_p^d$ is in one to one correspondence with the set of series
$$
x=\sum_{i=0}^{\infty}A^{i}x_i,
$$
where each $x_j$ runs over all possible digits from $\mathbb{Z}_p^d/A\mathbb{Z}_p^d$.

Analogously, the group
$\mathbb{Q}_p^d$ is in one to one correspondence with the set of series
\begin{equation}\label{expansion}
x=\sum_{i=\gamma}^{\infty}A^{i}x_i,\qquad \gamma\in \mathbb{Z},\quad x_i\in \mathbb{Z}_p^d/A\mathbb{Z}_p^d.
\end{equation}
}
\end{lemma}

\noindent{\it Proof}\qquad
The above series can be constructed taking residues of $x\in \mathbb{Q}_p^d$ with respect to the subgroups $A^{i}\mathbb{Z}_p^d$. Namely for $x\in \mathbb{Z}_p^d$ one has
$$
x_0=x\,{\rm mod}\, A\mathbb{Z}_p^d,\quad Ax_1=(x-x_0)\,{\rm mod}\, A^2\mathbb{Z}_p^d,\quad
A^2x_2=(x-x_0-Ax_1)\,{\rm mod}\, A^3\mathbb{Z}_p^d,\quad\dots.
$$
The existence and uniqueness (in the fixed set of representatives of equivalence classes) of solutions of the above equations for $x_i$ follows from the condition that $A$ is a dilation.

Convergence of the corresponding series (\ref{expansion}) is straightforward (since $A$ is a $\widetilde{d}$--contraction).

Ultrametricity of the deformed metric implies that different series correspond to different $x\in \mathbb{Q}_p^d$. If
$$
x=\sum_{i=\gamma}^{\infty}A^{i}x_i,\qquad x'=\sum_{i=\gamma}^{\infty}A^{i}x'_i
$$
and $i_0$ is the number of the first non--coinciding digit in the above expansions, i.e. $x_i=x'_i$, $i<i_0$, $x_{i_0}\ne x'_{i_0}$, then
by ultrametricity
$$
\|x-x'\|= \|A^{i_0}(x_{i_0}-x'_{i_0})\|.
$$
This finishes the proof of the lemma. $\square$

\bigskip

Let us note that the above parametrization is valid for any choice of the digits. Namely any choice of digits in the equivalence classes of $\mathbb{Z}_p^d/A\mathbb{Z}_p^d$ give the same sets $\mathbb{Z}_p^d$ and $\mathbb{Q}_p^d$ as the sets of the series of the above form.

In the following for the dilation matrix $A$ above we will claim that this matrix has integer matrix elements.

Let us consider the map of $\mathbb{Q}_p^d$ to some subset of $\mathbb{R}_p^d$ of the form
\begin{equation}\label{themap}
\rho:\sum_{i=\gamma}^{\infty}A^i x_i\mapsto
 \sum_{i=\gamma}^{\infty}A^{-i-1}x_i.
\end{equation}

This map and the image of $\mathbb{Z}_p^d$ with respect to this map do depend on the choice of the digits $\{n_0,\dots,n_{q-1}\}$. Since the digits and matrix elements of $A$ are (rational) integers the map is defined correctly.

In \cite{GroMad}, Theorem 3 the following equivalent conditions for the dilation $A$ (with the integer matrix elements) to be related to some real Haar multiresolution analysis were formulated:

(A)  The following subset of $\mathbb{R}^d$ given by the set of series
$$
R=\{\sum_{i=0}^{\infty}A^{-i-1}x_i\},
$$
where $x_i$ are the digits (integer representatives from $\mathbb{Z}^d/A\mathbb{Z}^d$), is Lebesgue measurable with the measure equal to one.

(B) The Lebesgue measure of the intersections $R\bigcap (R+k)$ for all non-zero $k\in\mathbb{Z}^d$ is equal to zero.

\begin{lemma}\quad
{\sl Let $A$ be a dilation (with the integer matrix elements) in $\mathbb{Q}_p^d$ with respect to the (deformed) metric $\widetilde{d}$. Let the corresponding map (\ref{themap}) conserve the measure (i.e. the Lebesgue measure of the image of a set is equal to the Haar measure of this set).

Then:

1) The above conditions (A), (B) of \cite{GroMad} are satisfied.

2) One has
$$
|{\rm det}\,A|_p=|{\rm det}\,A|^{-1}.
$$

}
\end{lemma}

\noindent{\it Proof}\qquad
We have $R=\rho(\mathbb{Z}_p^d)$.
If $\rho$ conserves the measure then the Haar measure $\mu(\mathbb{Z}_p^d)$ and the Lebesgue measure $|R|$ are equal to one. The statement $|R|=1$ is exactly the condition (A).

One has
\begin{equation}\label{rhoAnZ}
\rho(A^{j}(n+\mathbb{Z}_p^d))=A^{-j}\rho(n)+A^{-j}\rho(\mathbb{Z}_p^d)=A^{-j}\sum_{i=\gamma}^{-1}A^{-i-1} n_i+A^{-j}R.
\end{equation}
Here $j\in{\mathbb{Z}}$, $n\in \mathbb{Q}_p^d/\mathbb{Z}_p^d$, $n=\sum_{i=\gamma}^{-1}A^i n_i$ (where $n_i$ are the digits from $\mathbb{Z}_p^d/A\mathbb{Z}_p^d$)

Taking the Haar and the Lebesgue measures of $A^{j}(n+\mathbb{Z}_p^d)$ and $\rho(A^{j}(n+\mathbb{Z}_p^d))$ given by (\ref{rhoAnZ}) we get the necessary condition: if $\rho$ is a measure conserving map then
$$
|{\rm det}\,A|_p=|{\rm det}\,A|^{-1}.
$$

Namely the Lebesgue measure of (\ref{rhoAnZ}) is equal to $|{\rm det}\,A|^{-j}$.

If any pair of the sets $\rho(n+\mathbb{Z}_p^d)$, $\rho(n'+\mathbb{Z}_p^d)$, $n\ne n'$ given by (\ref{rhoAnZ}) possesses the intersection of a non-zero measure then the map $\rho$ can not be measure conserving.

This can be checked as follows: let us fix the minimal $\widetilde{d}$-ball in $\mathbb{Q}_p^d$ containing the both $n$ and $n'$. By lemma \ref{orbit} this ball is equal to some $A^{j}(m+\mathbb{Z}_p^d)$, $m\in \mathbb{Q}_p^d/\mathbb{Z}_p^d$, $j$ is negative. The ball $A^{j}(m+\mathbb{Z}_p^d)$ is a finite disjoint union of some translations of $\mathbb{Z}_p^d$ (the number of these subballs is equal to $|{\rm det}\,A|_p^{j}=|{\rm det}\,A|^{-j}$). The $\rho$-image of $A^{j}(m+\mathbb{Z}_p^d)$ is the (not necessarily disjoint) union of the the $\rho$-images of these translations of $\mathbb{Z}_p^d$. The Lebesgue measure $|\rho(A^{j}(m+\mathbb{Z}_p^d))|=|{\rm det}\,A|^{-j}$, the images of the translations of $\mathbb{Z}_p^d$ have the measure one.

If $\rho(n+\mathbb{Z}_p^d)$, $\rho(n'+\mathbb{Z}_p^d)$, $n\ne n'$ possess intersection of the non-zero measure, the sum of measures of $\rho$-images of the translations of $\mathbb{Z}_p^d$ will be less than $|{\rm det}\,A|^{-j}$ and the map $\rho$ will not conserve the measure.

Since both $A$ and $n_i$ have integer matrix elements, the term $A^{-j}\rho(n)$ for non-positive $j$ in (\ref{rhoAnZ}) belongs to $\mathbb{Z}^d$.

This implies some weaker form of the condition (B): the Lebesgue measure of the intersections $R\bigcap (R+k)$ for all non-zero $k\in\mathbb{Z}^d$, $k=\rho(n)$, $n\in \mathbb{Q}_p^d/\mathbb{Z}_p^d$ is equal to zero. Note that since by \cite{GroMad} (B) is equivalent to (A) (which we proved) the condition (B) holds for general $k\in\mathbb{Z}^d$.

This finishes the proof of the lemma. $\square$

\bigskip

Let us note that the condition $|{\rm det}\,A|_p=|{\rm det}\,A|^{-1}$ is not sufficient to make the map $\rho$ measure conserving --- the corresponding one-dimensional example was considered in the beginning of this section.

\begin{lemma}\quad{\sl
Let $A$ be a dilation (with the integer matrix elements) in $\mathbb{Q}_p^d$ with respect to the (deformed) metric $\widetilde{d}$.
Let the conditions (A), (B) of \cite{GroMad} be satisfied.

Then the map (\ref{themap}) conserves the measure (i.e. the Lebesgue measure of the image of a set is equal to the Haar measure of this set).
}
\end{lemma}

\noindent{\it Proof}\qquad
By (A) the image of $\mathbb{Z}_p^d$ has the measure one. By  (\ref{rhoAnZ})
$$
\mu(A^{j}(n+\mathbb{Z}_p^d)=|\rho(A^{j}(n+\mathbb{Z}_p^d)|=|{\rm det}\,A|_p^{j}=|{\rm det}\,A|^{-j},
$$
i.e. the measure of the image of any $\widetilde{d}$-ball in $\mathbb{Q}_p^d$ is equal to the Haar measure of this ball.

To prove that $\rho$ conserves the measure it is sufficient to prove that the Lebesgue measure of the image of a union of balls is equal to the Haar measure of this union of balls. This reduces to the following statement: the intersection of images of any two non-intersecting balls has the zero Lebesgue measure. It is sufficient to prove this for any pair of balls of the diameter one.

In the proof of the previous lemma it was found that this condition (zero measure intersection of the images of two balls of the diameter one) is guaranteed by the condition (B).
This finishes the proof of the lemma. $\square$

\bigskip

\noindent{\bf Acknowledgments}\qquad One of the authors (S.K.) would
like to thank I.V.Volovich for fruitful discussions. He gratefully
acknowledges being partially supported by the grant DFG Project 436
RUS 113/809/0-1 and DFG Project 436 RUS 113/951, by the grants of
the Russian Foundation for Basic Research
RFFI 08-01-00727-a and RFFI 09-01-12161-ofi-m, by the grant of the President of Russian
Federation for the support of scientific schools NSh-7675.2010.1 and
by the Program of the Department of Mathematics of the Russian
Academy of Science ''Modern problems of theoretical mathematics'',
and by the program of Ministry of Education and Science of Russia
''Development of the scientific potential of High School, years of
2009--2010'', project 3341. He is also grateful to IZKS (the
Interdisciplinary Center for Complex Systems) of the University of
Bonn for kind hospitality.

\end{document}